\documentclass[a4paper,11pt,reqno]{amsart}

\usepackage{amssymb}

\newtheorem{theorem}[equation]{Theorem}
\newtheorem{lemma}[equation]{Lemma}

\newtheorem{proposition}[equation]{Proposition}

\numberwithin{equation}{section}

\theoremstyle{definition}

\newtheorem*{example*}{Example}

\newtheorem{remark}[equation]{Remark}
\newtheorem*{remark*}{Remark}

\newcommand{\bZ}{{\mathbb Z}}

\newcommand{\cO}{{\mathcal O}}
\newcommand{\bj}{\textup{bj}}
\newcommand{\Bj}{\textup{Bj}}
\newcommand{\frg}{{\mathfrak g}}

\newcommand{\fre}{{\mathfrak e}}

\newcommand{\frs}{{\mathfrak s}}
\newcommand{\fra}{{\mathfrak a}}

\newcommand{\calT}{{\mathcal T}}

\newcommand{\calL}{{\mathcal L}}

\newcommand{\subo}{_{\bar 0}}
\newcommand{\subuno}{_{\bar 1}}

\providecommand{\espan}[1]{\text{span}\left\{ #1\right\}}

 \newcommand{\frsl}{{\mathfrak{sl}}}

 \newcommand{\frpsl}{{\mathfrak{psl}}}
 \newcommand{\frgl}{{\mathfrak{gl}}}

 \newcommand{\frd}{{\mathfrak{d}}}

 \DeclareMathOperator{\ad}{ad}
 \newcommand{\der}{\mathfrak{der}}
 \newcommand{\inder}{\mathfrak{inder}}
 \DeclareMathOperator{\End}{End}
 \DeclareMathOperator{\Hom}{Hom}
 \DeclareMathOperator{\Mat}{Mat}

\def\bigstrut{\vrule height 12pt width 0ptdepth 2pt}

\def\hreglon{\hrule height1pt}
\def\vreglon{\vrule height 12pt width1pt depth 4pt}

\def\hreglonfill{\leaders\hreglon\hfill}

\newenvironment{romanenumerate}
 {\begin{enumerate}
 
 }{\end{enumerate}}

\begin{document}

\title[Tits construction in characteristic $3$]{The Tits
construction and some simple Lie superalgebras in characteristic
$3$}

\author[Alberto Elduque]{Alberto Elduque$^{\star}$}
 \thanks{$^{\star}$ Supported by the Spanish Ministerio de
 Educaci\'{o}n y Ciencia
 and FEDER (MTM 2004-081159-C04-02) and by the
Diputaci\'on General de Arag\'on (Grupo de Investigaci\'on de
\'Algebra)}
 \address{Departamento de Matem\'aticas e Instituto Universitario de Matem\'{a}ticas y Aplicaciones,
  Universidad de Zaragoza, 50009 Zaragoza, Spain}
 \email{elduque@unizar.es}


\date{March 27, 2007}



\begin{abstract}
Some simple Lie superalgebras, specific of characteristic $3$,
defined by S.~Bouarroudj and D.~Leites \cite{BL}, will be related to
the simple alternative and commutative superalgebras discovered by
I.P.~Shestakov \cite{Shestakov}.
\end{abstract}

\maketitle



Throughout the paper, the ground field $k$ will always be assumed to
be of characteristic $\ne 2$.

\section{Tits construction}

In 1966 \cite{Tits66}, Tits gave a unified construction of the
exceptional simple classical Lie algebras by means of two
ingredients: a unital composition algebra and a degree three simple
Jordan algebra. The approach used by Benkart and Zelmanov in
\cite{BZ96} will be followed here (see also \cite{EOS4Tits}) to
review this construction.

Let $C$ be a unital composition algebra over the ground field $k$
with norm $n$. Thus, $C$ is a finite dimensional unital $k$-algebra,
with the nondegenerate quadratic form $n:C\rightarrow k$ such that
$n(ab)=n(a)n(b)$ for any $a,b\in C$. Then, each element satisfies
the degree $2$ equation
\begin{equation}\label{eq:deg2}
a^2-t(a)a+n(a)1=0,
\end{equation}
where $t(a)=n(a,1)\,\bigl(=n(a+1)-n(a)-1\bigr)$ is called the
\emph{trace}. The subspace of trace zero elements will be denoted by
$C^0$.

Moreover, for any $a,b\in C$, the linear map $D_{a,b}:C\rightarrow
C$ given by
\begin{equation}\label{eq:Dab}
D_{a,b}(c)=[[a,b],c]+3(a,c,b)
\end{equation}
where $[a,b]=ab-ba$ is the commutator, and $(a,c,b)=(ac)b-a(cb)$ the
associator, is a derivation: the \emph{inner derivation} determined
by the elements $a,b$ (see \cite[Chapter III, \S 8]{Schafer}). These
derivations satisfy
\begin{equation}\label{eq:Dcyclic}
D_{a,b}=-D_{b,a},\quad D_{ab,c}+D_{bc,a}+D_{ca,b}=0,
\end{equation}
for any $a,b,c\in C$. The linear span of these derivations will be
denoted by $\inder\, C$. It is an ideal of the whole Lie algebra of
derivations $\der\, C$ and, if the characteristic is $\ne 3$, it is
the whole $\der\, C$.

The dimension of $C$ is restricted to $1$, $2$, $4$ (quaternion
algebras) and $8$ (Cayley algebras). If the ground field $k$ is
algebraically closed, the only unital composition algebra are, up to
isomorphism, the ground field $k$, the cartesian product of two
copies of the ground field $k\times k$, the algebra of two by two
matrices $\Mat_2(k)$, and the split Cayley algebra $C(k)$. (See, for
instance, \cite[Chapter 2]{ZSSS}.)

\smallskip

Now, let $J$ be a unital Jordan algebra with a \emph{normalized
trace} $t_J:J\rightarrow k$. That is, $t_J$ is a linear map such
that $t_J(1)=1$ and $t_J\bigl((xy)z\bigr)=t_J\bigl(x(yz)\bigr)$ for
any $x,y,z\in J$. Then $J=k1\oplus J^0$, where $J^0=\{x\in J:
t_J(x)=0\}$. For any $x,y\in J^0$, the product $xy$ splits as
\begin{equation}\label{eq:Jxy}
xy=t_J(xy)1+x*y,
\end{equation}
with $x*y\in J^0$. Then $x*y=xy-t_J(xy)1$ gives a commutative
multiplication on $J^0$. The linear map $d_{x,y}:J\rightarrow J$
defined by
\begin{equation}\label{eq:dxy}
d_{x,y}(z)=x(yz)-y(xz),
\end{equation}
is the inner derivation of $J$ determined by the elements $x$ and
$y$. Since $d_{1,x}=0$ for any $x$, it is enough to deal with the
inner derivations $d_{x,y}$, with $x,y\in J^0$. The linear span of
these derivations will be denoted by $\inder\, J$, which is an ideal
of the whole Lie algebra of derivations $\der\, J$.

\smallskip

Given $C$ and $J$ as before, consider the space
\begin{equation}\label{eq:TCJ}
\calT(C,J)=\inder\, C\oplus \bigl(C^0\otimes J^0\bigr)\oplus
\inder\, J
\end{equation}
(unadorned tensor products are always considered over $k$), with the
anticommutative multiplication $[.,.]$ specified by:
\begin{equation}\label{eq:TCJproduct}
\begin{split}
\bullet&\ \textrm{$\inder\, C$ and $\inder\, J$ are Lie subalgebras,}\\
\bullet&\ [\inder\, C,\inder\, J]=0,\\
\bullet&\ [D,a\otimes x]=D(a)\otimes x,\ [d,a\otimes x]=a\otimes
d(x),\\
\bullet&\ [a\otimes x,b\otimes y]=t_J(xy)D_{a,b}+\bigl([a,b]\otimes
x*y\bigr)+2t(ab)d_{x,y},
\end{split}
\end{equation}
for all $D\in \inder\, C$, $d\in \inder\, J$, $a,b\in C^0$, and
$x,y\in J^0$.

\smallskip

The conditions for $\calT(C,J)$ to be a Lie algebra are the
following:
\begin{equation}\label{eq:TCJLie}
\begin{split}
\textrm{(i)}&\  \displaystyle{\sum_{\circlearrowleft}
 t\bigl([a_{1}, a_{2}] a_{3}\bigr)\,
d_{(x_1 * x_2), x_3}}=0,\\[6pt]
\textrm{(ii)}&\  \displaystyle{\sum_{\circlearrowleft}
 t_J\bigl( (x_1 * x_2) x_{3}\bigr)
\,D_{[a_1, a_2], a_3}}=0,\\[6pt]
\textrm{(iii)}&\ \displaystyle{\sum_{\circlearrowleft}
 \Bigl(D_{a_1,a_2}(a_3) \otimes t_J\bigl(x_1
x_2\bigr) x_3} + [[a_1, a_2],a_3] \otimes (x_1 * x_2)* x_3\\[-6pt]
  &\qquad\qquad\qquad\qquad +2
t(a_1 a_2) a_3\otimes d_{x_1, x_2}(x_3)\Bigr)=0
\end{split}
\end{equation}
for any $a_1,a_2,a_{3} \in C^{0}$ and any $x_1,x_2,x_3 \in J^0$. The
notation ``$\displaystyle{\sum_\circlearrowleft}$'' indicates
summation over the cyclic permutation of the variables.

These conditions appear in \cite[Proposition 1.5]{BE03}, but there
they are stated in the more general setting of superalgebras, a
setting we will deal with later on. In particular, over fields of
characteristic $\ne 3$, these conditions are fulfilled if $J$ is a
separable Jordan algebra of degree three over $k$ and
$t_J=\frac{1}{3}T$, where $T$ denotes the generic trace of $J$ (see
for instance \cite{Jac68}).

Over an algebraically closed field $k$ of characteristic $\ne 3$,
the degree $3$ simple Jordan algebras are, up to isomorphism, the
algebras of $3\times 3$ hermitian matrices over a unital composition
algebra: $H_3(C')$ (see \cite{Jac68}). By varying $C$ and $C'$,
$\calT(C,H_3(C'))$ is a classical simple Lie algebra, and
Freudenthal's Magic Square (Table \ref{ta:FMS}) is obtained.

\begin{table}[h!]
$$
\vbox{\offinterlineskip
 \halign{\hfil\ $#$\ \hfil&%
 \vreglon #%
 &\hfil\ $#$\ \hfil&\hfil\ $#$\ \hfil
 &\hfil\ $#$\ \hfil&\hfil\ $#$\ \hfil\cr
 \bigstrut &width 0pt&H_3(k)&H_3(k\times k)&H_3(\Mat_2(k))&H_3(C(k))\cr
 &\multispan5{\hreglonfill}\cr
 k&&A_1&A_2&C_3&F_4\cr
 \bigstrut k\times k&& A_2&A_2\oplus A_2&A_5&E_6\cr
 \bigstrut \Mat_2(k)&&C_3 & A_5&D_6&E_7\cr
 \bigstrut C(k)&& F_4& E_6& E_7&E_8\cr}}
$$
\bigskip
\caption{Freudenthal's Magic Square}\label{ta:FMS}
\end{table}

\bigskip

Let us have a look at the rows in the Tits construction of
Freudenthal's Magic Square.

\medskip

\noindent\textbf{First row:}\quad Here $C=k$, so $C^0=0$ and
$\inder\, C=0$. Thus, $\calT(C,J)$ is just $\inder\, J$. In
particular, $\calT(k,J)$ makes sense and is a Lie algebra for any
Jordan algebra $J$.

\medskip

\noindent\textbf{Second row:}\quad Here $C=k\times k$, so $C^0$
consists of the scalar multiples of $(1,-1)$, and thus $\calT(C,J)$
can be identified with $J^0\oplus \inder\, J$. The elements in $J^0$
multiply as $[x,y]=4d_{x,y}$ because
$t\bigl((1,-1)^2\bigr)=t\bigl((1,1)\bigr)=2$. Given any Jordan
algebra $J$, its \emph{Lie multiplication algebra} $\calL(J)$ (see
\cite{Schafer}) is the Lie subalgebra of the general linear Lie
algebra $\frgl(J)$ generated by $l_J=\{l_x:x\in J\}$, where
$l_x:y\mapsto xy$ denotes the left multiplication by $x$. Then
$\calL(J)=l_J\oplus\inder\, J$. The map
\[
\begin{split}
\calT(C,J)&\rightarrow \calL(J)\\
(1,-1)\otimes x+d&\mapsto 2l_x+d,
\end{split}
\]
is a monomorphism. Its image is $\calL^0(J)=l_{J^0}\oplus \inder\,
J$. Again this shows that $\calT(k\times k,J)$ makes sense and is a
Lie algebra for any Jordan algebra with a normalized trace. Given
any separable Jordan algebra of degree $3$ over a field $k$ of
characteristic $\ne 3$, $\calL^0(J)$ is precisely the derived
algebra $[\calL(J),\calL(J)]$. This latter Lie algebra makes sense
for any Jordan algebra over any field. (Recall that the
characteristic is assumed to be $\ne 2$ throughout.)

\medskip

\noindent\textbf{Third row:}\quad Here $C=\Mat_2(k)$ or, if the
ground field is not assumed to be algebraically closed, $C$ is any
quaternion algebra $Q$. Under these circumstances, $Q^0$ is a simple
three-dimensional Lie algebra under the commutator ($[a,b]=ab-ba$),
and any simple three-dimensional Lie algebra appears in this way.
Besides, for any $a,b\in Q^0$, the inner derivation $D_{a,b}$ is
just $\ad_{[a,b]}$, since $Q$ is associative. Hence, $\inder\, Q$
can be identified with $Q^0$, and $\calT(Q,J)$ with
\[
Q^0\oplus \bigl(Q^0\otimes J^0\bigr)\oplus\inder\,
J\simeq\bigl(Q^0\otimes (k1\oplus J^0)\bigr)\oplus\inder\, J\simeq
\bigl(Q^0\otimes J\bigr)\oplus \inder\, J,
\]
and the Lie bracket \eqref{eq:TCJproduct} in $\calT(Q,J)$ becomes
the bracket in $\bigl(Q^0\otimes J\bigr)\oplus \inder\, J$ given by
\[
\begin{split}
\bullet&\ \textrm{$\inder\, J$ is a Lie subalgebra,}\\
\bullet&\ [d,a\otimes x]=a\otimes d(x),\\
\bullet&\ [a\otimes x,b\otimes y]=\bigl([a,b]\otimes
xy\bigr)+2t(ab)d_{x,y},
\end{split}
\]
for any $a,b\in Q^0$, $x,y\in J$, and $d\in\inder\, J$, since
$t_J(xy)1+x*y=xy$ for any $x,y\in J$. This bracket makes sense for
any Jordan algebra (not necessarily endowed with a normalized
trace), it goes back to \cite{Tits62} and, in the split case
$Q=\Mat_2(k)$, the resulting Lie algebra is the well-known
Tits-Kantor-Koecher Lie algebra of the Jordan algebra $J$.

\medskip

\noindent\textbf{Fourth row:}\quad In the last row, $C$ is a Cayley
algebra over $k$. If the characteristic of the ground field $k$ is
$\ne 3$, the Lie algebra $\der\, C=\inder\, C$ is a simple Lie
algebra of type $G_2$ (dimension $14$), and $C^0$ is its unique
seven dimensional irreducible module. In particular, over any
algebraically closed field of characteristic $\ne 3$, the Lie
algebra $\calT\bigl(C(k),J\bigr)$ is a Lie algebra graded over the
root system $G_2$. These $G_2$-graded Lie algebras contain a simple
subalgebra isomorphic to $\der\, C(k)$ such that, as modules for
this subalgebra, they are direct sums of copies of modules of three
types: adjoint, the irreducible seven dimensional module, and the
trivial one dimensional module. These Lie algebras have been
determined in \cite{BZ96} and the possible Jordan algebras involved
are essentially the degree $3$ Jordan algebras.

In characteristic $3$, however, the situation is completely
different. To begin with, given a Cayley algebra $C$ over a field
$k$ of characteristic $3$, its Lie algebra of derivations $\der\, C$
is no longer simple (see \cite{AME}) but contains a unique minimal
ideal, which is precisely $\inder\, C=\ad_{C^0}$ (note that
$D_{a,b}=\ad_{[a,b]}$ in characteristic $3$ because of
\eqref{eq:Dab}), which is isomorphic to the Lie(!) algebra
$\bigl(C^0,[.,.]\bigr)$. This latter Lie algebra is a form of the
projective special linear Lie algebra $\frpsl_3(k)$ (and any form of
$\frpsl_3(k)$ appears in this way \cite[\S 4]{Egroctonions}).
Moreover, in \cite{AME} it is shown that the quotient $\der\,
C/\inder\, C=\der\, C/\ad_{C^0}$ is isomorphic too, as a Lie
algebra, to $\bigl(C^0,[.,.]\bigr)$. Therefore, the algebra
$\calT(C,J)$ in \eqref{eq:TCJ} can be identified in this case with
\[
\begin{split}
C^0\oplus\bigl(C^0\otimes J^0\bigr)\oplus\inder\, J&\simeq
 \bigl(C^0\otimes (k1\oplus J)\bigr)\oplus\inder\, J\\
 &\simeq(C^0\otimes J)\oplus \inder\, J,
\end{split}
\]
and hence, as a module for $\inder\, C\simeq C^0$, it is a direct
sum of copies of the adjoint module and of the trivial module. The
Lie algebras and superalgebras with these properties  will be
determined in this paper.

\section{A family of Lie algebras}

Throughout this section, the ground field $k$ will be always assumed
to be of characteristic $3$. Let $C$ be a Cayley algebra, that is,
an eight dimensional unital composition algebra over $k$, and let
$C^0$ denote its subspace of trace zero elements.

For any $a,b,c\in c$, a simple computation using that the associator
is skew symmetric on its arguments since $C$ is alternative gives:
\[
\begin{split}
[[a,b],c&+[[b,c],a]+[[c,a],b]\\
 &=(a,b,c)-(b,a,c)+(b,c,a)-(c,b,a)+(c,a,b)-(a,c,b)\\
 &=6(a,b,c)=0.
\end{split}\]
Hence $C$ is a Lie algebra under the bracket $[a,b]=ab-ba$, and
$C^0$ is an ideal of $C$. Besides, for any $a,b\in C^0$, since the
subalgebra generated by any two elements in $C$ is associative
(Artin's Theorem, see \cite{Schafer}), one obtains
\[
\begin{aligned}
\null [[a,b],b]&=ab^2+b^2a-2bab&&\\
 &=ab^2+b(ba+ab)&&\text{(as $2=-1$)}\\
 &=-n(b)a-n(a,b)b&&\text{(as $a^2=-n(a)1$ for any $a\in C^0$)}
\end{aligned}
\]
Thus, for any $a,b\in C^0$,
\begin{equation}\label{eq:[abb]}
[[a,b],b]=n(b,b)a-n(a,b)b.
\end{equation}

If $C$ is split, then it contains a basis
$\{e_1,e_2,u_1,u_2,u_3,v_1,v_2,v_3\}$ with multiplication given by
(see \cite{Schafer,ZSSS}):
\begin{equation}\label{eq:productC}
\begin{split}
&e_i^2=e_i,\ i=1,2,\quad e_1e_2=e_2e_1=0,\\
&e_1u_j=u_j=u_je_2,\ e_2v_j=v_j=v_je_1,\ j=1,2,3,\\
&e_2u_j=u_je_1=0=e_1v_j=v_je_2,\ j=1,2,3,\\
&u_iv_j=-\delta_{ij}e_1,\ v_iu_j=-\delta_{ij}e_2,\ i,j=1,2,3\
\text{($\delta_{ij}$ is $1$ for $i=j$, $0$ otherwise)}\\
&u_iu_j=\epsilon_{ijk}v_k,\ v_iv_j=\epsilon_{ijk}u_k,\
\text{($\epsilon_{ijk}$ skew symmetric with $\epsilon_{123}=1$).}
\end{split}
\end{equation}
Moreover,
\begin{equation}\label{eq:normC}
\begin{split}
&n(e_i)=0=n(u_j,u_k)=n(v_j,v_k),\ i=1,2,\, j,k=1,2,3,\\
&n(e_1,e_2)=1,\ n(u_j,v_k)\delta_{jk},\ j,k=1,2,3.
\end{split}
\end{equation}
Then, with $h=e_1-e_2$, $C^0$ is the linear span of
$\{h,u_1,u_2,u_3,v_1,v_2,v_3\}$, and these elements multiply as:
\begin{equation}\label{eq:bracketC0}
\begin{split}
&[h,u_j]=2u_j,\ [h,v_j]=-2v_j,\ j=1,2,3,\\
&[u_i,u_j]=2\epsilon_{ijk}v_k,\ [v_i,v_j]=2\epsilon_{ijk}u_k,\
i,j,k=1,2,3,\\
&[u_j,v_k]=-\delta_{jk}h,\ j,k=1,2,3.
\end{split}
\end{equation}

Denote by $\frs$ the Lie algebra $\bigl(C^0,[.,.]\bigr)$. It is easy
to check that the $\frs$-module $\frs\otimes\frs$ is generated by
$u_1\otimes v_1$, which is an eigenvector for $\ad_h$ with
eigenvalue $0$. From this fact, it follows at once that the
dimension of the space of invariant linear maps
$\Hom_{\frs}(\frs\otimes\frs,\frs)$ is $1$, being this space
generated by the Lie bracket. Also, $\Hom_{\frs}(\frs\otimes\frs,k)$
is one dimensional, spanned by the bilinear map induced by the norm
$n(.,.)$. By extending scalars, this is shown to be valid for any
Cayley algebra, not necessarily split.

\smallskip

Thus, let $C$ be any Cayley algebra over $k$ and let $\frs$ be the
Lie algebra $\bigl(C^0,[,.,]\bigr)$. Let $\frg$ be a Lie algebra
endowed with an action of $\frs$ on $\frg$ by derivations:
$\rho:\frs\rightarrow \der\, \frg$, such that, as a module for
$\frs$, $\frg$ is a direct sum of copies of the adjoint module and
the one dimensional trivial module. Gathering together the copies of
the adjoint module and the copies of the trivial module, $\frg$ can
be identified with
\begin{equation}\label{eq:gsAd}
\frg=\bigl(\frs\otimes A\bigr)\oplus \frd,
\end{equation}
where $\frd$ is the sum of the trivial $\frs$-modules and $A$ is a
vector space. As $\frd$ equals $\{d\in \frg : \rho(s)(d)=0\ \forall
s\in \frs\}$, it is a subalgebra of $\frs$. Now, the invariance of
the bracket in $\frg$ under the action of $\frs$, together with the
fact that $\Hom_{\frs}(\frs\otimes\frs,\frs)$ (respectively
$\Hom_{\frs}(\frs\otimes\frs,k)$) is spanned by the Lie bracket
(respectively, the bilinear form induced by the norm) shows that
there are bilinear maps
\begin{equation}\label{eq:bilinear}
\begin{aligned}
\frd\times A&\rightarrow A,& (d,a)&\mapsto d(a),\\
A\times A&\rightarrow A,& (a_1,a_2)&\mapsto a_1a_2,\\
A\times A&\rightarrow \frd,& (a_1,a_2)&\mapsto d_{a_1,a_2},
\end{aligned}
\end{equation}
such that the Lie bracket on $\frg=\bigl(\frs\otimes
A\bigr)\oplus\frd$ is given by:
\begin{equation}\label{eq:product}
\begin{split}
\bullet&\ [d,s\otimes a]=s\otimes d(a)=-[s\otimes a,d],\\
\bullet&\ [s_1\otimes a_1,s_2\otimes a_2]=[s_1,s_2]\otimes
a_1a_2+n(s_1,s_2)d_{a_1,a_2},\\
\bullet&\ [d_1,d_2]\ \text{is the product in the subalgebra $\frd$},
\end{split}
\end{equation}
for any $d,d_1,d_2\in \frd$, $s,s_1,s_2\in \frs$, and $a,a_1,a_2\in
A$. The skew symmetry of the Lie bracket forces the product $a_1a_2$
on $A$ to be commutative, and the bilinear map $d_{a_1,a_2}$ to be
skew symmetric.

\smallskip

Now, let us consider the Jacobi identity $J(z_1,z_2,z_3)=0$ on
$\frg$, where $J(z_1,z_2,z_3)=\displaystyle{\sum_{\circlearrowleft}
 [[z_{1}, z_{2}] z_{3}]}$:

\smallskip

\begin{itemize}

\item With $z_1=d_1$, $z_2=d_2$ in $\frd$ and $z_3=s\otimes a$,
$s\in \frs$, $a\in A$, this gives
$[d_1,d_2](a)=d_1(d_2(a))-d_2(d_1(a))$. That is, the linear map
$\Phi:\frd\rightarrow \frgl(A)$, $\Phi(d):a\mapsto d(a)$, is a
representation of the Lie algebra $\frd$.

\smallskip

\item With $z_1=d$, $z_2=s_1\otimes a_1$ and $z_3=s_2\otimes a_2$,
$d\in\frd$, $s_1,s_2\in\frs$ and $a_1,a_2\in A$, the Jacobi identity
gives:
\[
\begin{split}
&d(a_1a_2)=d(a_1)a_2+a_1d(a_2),\\
&[d,d_{a_1,a_2}]=d_{d(a_1),a_2}+d_{a_1,d(a_2)},
\end{split}
\]
for any $a_1,a_2\in A$. That is, $\Phi(\frd)\subseteq \der\, A$
holds and $d:A\times A\rightarrow \frd$ is a $\frd$-invariant
bilinear map.

\smallskip

\item With $z_i=s_i\otimes a_i$, $i=1,2,3$, $s_i\in \frs$, $a_i\in
A$, the Jacobi identity gives:
\begin{subequations}\label{eq:sotimesAd}
\begin{align}
&\displaystyle{\sum_{\circlearrowleft}
    n\bigl([s_1,s_2],s_3\bigr)d_{a_1a_2,a_3}=0}\label{eq:d}\\
&\displaystyle{\sum_{\circlearrowleft}\Bigl(
\bigl([[s_1,s_2],s_3]\otimes (a_1a_2)a_3\bigr) +
\bigl(n(s_1,s_2)s_3\otimes
d_{a_1,a_2}(a_3)\bigr)\Bigr)=0}\label{eq:sotimesA}
\end{align}
\end{subequations}

Extending scalars, it can be assumed that $C$ is split, so that a
basis of $C^0$ as in \eqref{eq:bracketC0} is available. With
$s_i=u_i$, $i=1,2,3$, $n([u_1,u_2],u_3)=2n(v_3,u_3)=2$ and
cyclically. (Note that, because of the invariance under derivations
of the norm, $n([s_1,s_2],s_3)=n(s_1,[s_2,s_3])=n([s_2,s_3],s_1)$.)
Hence \eqref{eq:d} is equivalent to
\begin{equation}\label{eq:dcyclic}
d_{a_1a_2,a_3}+d_{a_2a_3,a_1}+d_{a_3a_1,a_2}=0
\end{equation}
for any $a_1,a_2,a_3\in A$.

Also, $[[u_1,u_2],u_3]=2[v_3,u_3]=2h$ and cyclically, and
$n(u_j,u_k)=0$ for any $j,k$. Thus \eqref{eq:sotimesA} gives
\begin{equation}\label{eq:jacobiA}
(a_1a_2)a_3+(a_2a_3)a_1+(a_3a_1)a_2=0
\end{equation}
for any $a_1,a_2,a_3\in A$.

With $s_1=u_1$, $s_2=v_1$ and $s_3=h$, $[[s_1,s_2],s_3]=0$,
$[[s_2,s_3],s_1]=2h$, and $[[s_3,s_1],s_2]=-2h$, while
$n(s_1,s_2)=1$, and $n(s_2,s_3)=0=n(s_3,s_1)$, so
\eqref{eq:sotimesA} gives (note that $2=-1$ in $k$):
\[
(a_3a_1)a_2-(a_2a_3)a_1+d_{a_1,a_2}(a_3)=0,
\]
which, by the commutativity of the product on $A$, is equivalent to:
\begin{equation}\label{eq:da1a2}
d_{a_1,a_2}(a)=a_1(a_2a)-a_2(a_1a),
\end{equation}
for any $a_1,a_2,a\in A$.

\end{itemize}

\medskip

\begin{lemma}\label{le:jacobialternative}
Let $k$ be a field of characteristic $3$. The commutative algebras
over $k$ satisfying \eqref{eq:jacobiA} are precisely the commutative
alternative algebras. Moreover, given any such algebra $A$, for any
$a_1,a_2$ consider the linear map $d_{a_1,a_2}=[l_{a_1},l_{a_2}]$,
where $l_a$ denotes the multiplication by $a$. Then $d_{a_1,a_2}$ is
a derivation of $A$ and equation \eqref{eq:dcyclic} is satisfied.
\end{lemma}
\begin{proof} By commutativity, \eqref{eq:jacobiA} is equivalent to
$2(a_1a_2)a_2 + a_1a_2^2=0$, or ($2=-1$) to $a_1a_2^2=(a_1a_2)a_2$
for any $a_1,a_2$, which is the right alternative law. Because of
the commutativity, this is equivalent to the left alternative law,
and hence the algebra is alternative. Now, any commutative
alternative algebra is a Jordan algebra, since the Jordan identity
is $(x^2,y,x)=0$ for any $x,y$, which is satisfied because any two
elements in an alternative algebra generate an associative
subalgebra (Artin's Theorem, see \cite{Schafer}). Hence
$d_{a_1,a_2}=[l_{a_1},l_{a_2}]$ is a derivation of $A$. Finally,
equation \eqref{eq:dcyclic} becomes the linearization of
$[l_{x^2},l_x]=0$.
\end{proof}

\bigskip

Conversely, let $C$ be a Cayley algebra over $k$ and let $\frs$ be
the Lie algebra $\bigl(C^0,[.,.]\bigr)$. Let $A$ be a commutative
alternative algebra, and let $\frd$ be a Lie algebra endowed with a
Lie algebra homomorphism $\Phi:\frd\rightarrow \der\, A$ (thus, in
particular, $A$ is a module for $\frd$), and a $\frd$-invariant
bilinear map $d:A\times A\rightarrow \frd$, $(a_1,a_2)\mapsto
d_{a_1,a_2}$, such that $\Phi(d_{a_1,a_2})=[l_{a_1},l_{a_2}]$ and
$\displaystyle{\sum_{\circlearrowleft} d_{a_1a_2,a_3}}=0$ for any
$a_1,a_2,a_3\in A$. Then equation \eqref{eq:d} holds trivially by
the invariance of the norm in $C$, and equation \eqref{eq:sotimesA}
holds too, as it is equivalent to
\[
\begin{split}
0&=\displaystyle{\sum_{\circlearrowleft}\Bigl(
[[s_1,s_2,],s_3]\otimes(a_1a_2)a_3
+n(s_1,s_2)s_3\otimes\bigl(a_1(a_2a_3)-a_2(a_1a_3)\bigr)\Bigr)}\\
&=\displaystyle{\sum_{\circlearrowleft}\Bigl(
\bigl([[s_1,s_2],s_3]-n(s_2,s_3)s_1+n(s_3,s_1)s_2\bigr)\otimes
(a_1a_2)a_3\Bigr)}.
\end{split}
\]
But $(a_1a_2)a_3=-(a_2a_3)a_1-(a_3a_1)a_2$, so \eqref{eq:sotimesAd}
holds if
\[
[[s_1,s_2],s_3]-n(s_2,s_3)s_1+n(s_3,s_1)s_2
=[[s_2,s_3],s_1]-n(s_3,s_1)s_2+n(s_1,s_2)s_3
\]
for any $s_1,s_2,s_3\in \frs$, or
\[
[[s_1,s_2],s_3]+[[s_2,s_3],s_1]=2n(s_1,s_3)s_2-n(s_2,s_3)s_1-n(s_2,s_1)s_3,
\]
which is equivalent to \eqref{eq:[abb]}.

Therefore:

\begin{theorem}\label{th:sAd}
Let $C$ be a Cayley algebra over a field $k$ of characteristic $3$,
let $\frs$ be the Lie algebra $\bigl(C^0,[.,.]\bigr)$. Let $\frg$ be
a Lie algebra with an action of $\frs$ by derivations such that, as
a module for $\frs$, $\frg$ is a direct sum of irreducible modules
of two types: the adjoint and the trivial one-dimensional modules.
Then there is a commutative alternative algebra $A$ over $k$ and a
Lie algebra $\frd$ over $k$, endowed with a Lie algebra homomorphism
$\Phi:\frd\rightarrow \der\, A$ and a $\frd$-invariant skewsymmetric
bilinear map $d:A\times A\rightarrow \frd$, $(a_1,a_2)\mapsto
d_{a_1,a_2}$ with $\Phi(d_{a_1,a_2})=[l_{a_1},l_{a_2}]$ for any
$a_1,a_2\in A$, such that $\frg$ is isomorphic to the Lie algebra
\[
\bigl(\frs\otimes A)\oplus\frd
\]
with Lie bracket given by
\begin{equation}\label{eq:sAdproduct}
\begin{split}
\bullet&\ \textrm{$\frd$ is a Lie subalgebra,}\\[2pt]
\bullet&\ [d,s\otimes a]=s\otimes d(a),\ \textrm{for $s\in\frs$, $a\in A$,
$d\in\frd$,}\\[2pt]
\bullet&\ [s_1\otimes a_1,s_2\otimes a_2]= [s_1,s_2]\otimes a_1a_2
+n(s_1,s_2)d_{a_1,a_2},\\
&\qquad\qquad \textrm{for $s_1,s_2\in\frs$, and $a_1,a_2\in A$,.}
\end{split}
\end{equation}

Conversely, the formulas in \eqref{eq:sAdproduct} define a Lie
algebra on the vector space $\bigl(\frs\otimes A\bigr)\oplus\frd$,
which is endowed with an action of $\frs$ by derivations:
$\rho:\frs\rightarrow \der\bigl((\frs\otimes A)\oplus\frd\bigr)$,
such that $\rho(s)(s'\otimes a)=[s,s']\otimes a$, $\rho(s)(d)=0$,
for any $s,s'\in S$, $a\in A$ and $d\in \frd$.
\end{theorem}

\bigskip

\begin{remark}
Over fields of characteristic $\ne 2,3$, any commutative alternative
algebra is associative, because for any $x,y,z$,
\[
\begin{split}
3(x,y,z)&=(x,y,z)+(y,z,x)+(z,xy)\\
&=(xy)z-x(yz)+(yz)x-y(zx)+(zx)y-z(xy)\\
&=[xy,z]+[yz,x]+[zx,y]=0.
\end{split}
\]
\end{remark}

\bigskip

\begin{remark}
Let $A$ be a unital commutative alternative algebra over a field $k$
of characteristic $3$ such that $1\not\in (A,A,A)$ (this is the case
for the unital commutative associative algebras). Then $A$ is a
Jordan algebra with a normalized trace, because if $A^0$ is any
codimension $1$ subspace of $A$ containing $(A,A,A)$ but not
containing $1$, then $A=k1\oplus A^0$, and the linear form
$t:A\rightarrow k$, such that $t(1)=1$ and $t(A^0)=0$ is a
normalized trace. Let $C$ be a Cayley algebra over $k$ and define
$\calT(C,A)$ as in \eqref{eq:TCJ}:
\[
\calT(C,A)=\inder\, C\oplus\bigl(C^0\otimes A^0\bigr)\oplus \inder\,
A.
\]
Then $\calT(C,A)$ is a Lie algebra (bracket as in
\eqref{eq:TCJproduct}, isomorphic to $\bigl(\frs\otimes
A\bigr)\oplus \frd$, with $\frd=\inder\, A=[l_A,l_A]$, $\frs=C^0$,
and bracket as in Theorem \ref{th:sAd} (with $\Phi$ the natural
inclusion).

This gives the natural extension of the fourth row in Tits
construction to characteristic $3$. The Jordan algebras that appear
have nothing to do with the separable degree $3$ Jordan algebras.

It must be remarked that the simple commutative alternative algebras
are just the fields, but there are nontrivial prime commutative
alternative algebras in characteristic $3$ (see \cite{ZSSS}). Recall
that an algebra is simple if its multiplication is not trivial and
it contains no proper ideal, while it is prime if the product of any
two nonzero ideals is again nonzero.
\end{remark}

\section{superalgebras}

All the arguments used in the proofs of Lemma
\ref{le:jacobialternative} and Theorem \ref{th:sAd} are valid in the
setting of superalgebras, if parity signs are added suitably. The
super version of Theorem \ref{th:sAd} is:

\begin{theorem}\label{th:supersAd}
Let $C$ be a Cayley algebra over a field $k$ of characteristic $3$,
let $\frs$ be the Lie algebra $\bigl(C^0,[.,.]\bigr)$. Let $\frg$ be
a Lie superalgebra with an action of $\frs$ by (even) derivations
such that, as a module for $\frs$, $\frg$ is a direct sum of
irreducible modules of two types: adjoint and trivial. Then there is
a commutative alternative superalgebra $A$ over $k$ and a Lie
superalgebra $\frd$ over $k$, endowed with a homomorphism of Lie
superalgebras $\Phi:\frd\rightarrow \der\, A$ and an even
$\frd$-invariant (relative to $\Phi$) super skewsymmetric bilinear
map $d:A\times A\rightarrow \frd$, $(a_1,a_2)\mapsto d_{a_1,a_2}$
with $\Phi(d_{a_1,a_2})=[l_{a_1},l_{a_2}]$ for any $a_1,a_2\in A$,
such that $\frg$ is isomorphic to the Lie superalgebra
\[
\bigl(\frs\otimes A)\oplus\frd
\]
with Lie bracket given by
\begin{equation}\label{eq:supersAdproduct}
\begin{split}
\bullet&\ \textrm{$\frd$ is a Lie subalgebra,}\\[2pt]
\bullet&\ [d,s\otimes a]=s\otimes d(a),\ \textrm{for $s\in\frs$, $a\in A$,
$d\in\frd$,}\\[2pt]
\bullet&\ [s_1\otimes a_1,s_2\otimes a_2]= [s_1,s_2]\otimes a_1a_2
+n(s_1,s_2)d_{a_1,a_2},\\
&\qquad\qquad \textrm{for $s_1,s_2\in\frs$, and $a_1,a_2\in A$.}
\end{split}
\end{equation}

Conversely, the formulas in \eqref{eq:supersAdproduct} define a Lie
superalgebra on the vector superspace $\bigl(\frs\otimes
A\bigr)\oplus\frd$ (the even part is $\bigl(\frs\otimes
A\subo\bigr)\oplus \frd\subo$ and the odd part is $\bigl(\frs\otimes
A\subuno\bigr)\oplus \frd\subuno$), which is endowed with an action
of $\frs$ by (even) derivations: $\rho:\frs\rightarrow
\der\bigl((\frs\otimes A)\oplus\frd\bigr)$, such that
$\rho(s)(s'\otimes a)=[s,s']\otimes a$, $\rho(s)(d)=0$, for any
$s,s'\in S$, $a\in A$ and $d\in \frd$.
\end{theorem}

\medskip

Recall that given a superalgebra $A$, its Lie superalgebra of
derivations is the Lie superalgebra $\der\, A=(\der\,
A)\subo\oplus(\der\, A)\subuno$ (a subalgebra of the general linear
Lie superalgebra $\frgl(A)$), where for any homogeneous $d\in\der\,
A$ and homogeneous $a_1,a_2\in A$:
\[
d(a_1a_2)=d(a_1)a_2+(-1)^{da_1}a_1d(a_2),
\]
where, as usual, $(-1)^{da_1}$ is $-1$ if both $d$ and $a_1$ are
odd, and $(-1)^{da_1}$ is $1$ otherwise. The Lie bracket of
homogeneous elements in $\frgl(A)$ is given by
$[f,g]=fg-(-1)^{fg}gf$. The fact that $d$ in the Theorem above is
\emph{even} means that $d(A_i,A_j)$ is contained in $\frd_{i+j}$ for
any $i,j\in\{\bar 0,\bar 1\}$; and the invariance of $d$ relative to
$\Phi$ means that
\[
[f,d_{a_1,a_2}]=d_{\Phi(f)(a_1),a_2}+(-1)^{da_1}d_{a_1,\Phi(f)(a_2)}
\]
for any homogeneous elements $f\in \frd$ and $a_1,a_2\in A$.

\bigskip

The importance of Theorem \ref{th:supersAd} lies in the fact that
there do exist interesting examples of commutative alternative
simple superalgebras in characteristic $3$. Besides:

\begin{proposition}
Let $A$ be a nonzero commutative alternative superalgebra over a
field $k$ of characteristic $3$, and let $\frd$ be a Lie
superalgebra endowed with a homomorphism of Lie superalgebras
$\Phi:\frd\rightarrow \der\, A$ and an invariant (relative to
$\Phi$) skewsymmetric bilinear map $d:A\times A\rightarrow \frd$
with $\Phi(d_{a_1,a_2})=[l_{a_1},l_{a_2}]$ for any $a_1,a_2\in A$.
Let $\frg=\bigl(\frs\otimes A\bigr)\oplus\frd$ be the Lie
superalgebra constructed by means of \eqref{eq:supersAdproduct}.
Then $\frg$ is simple if and only if the following conditions are
fulfilled:
\begin{romanenumerate}
\item $\Phi$ is one-to-one,
\item $\frd=d_{A,A}\, (=\espan{d_{a_1,a_2}: a_1,a_2\in A})$,
\item $A$ is simple.
\end{romanenumerate}
\end{proposition}
\begin{proof}
Assume first that $\frg$ is simple. Since $\ker\Phi$ is an ideal,
not only of $\frd$, but of the whole $\frg$, it follows that $\Phi$
is one-to-one. Also, $\bigl(\frs\otimes A\bigr)\oplus d_{A,A}$ is an
ideal of $\frg$, so $d_{A,A}=\frd$. Finally, if $I$ is a nonzero
ideal of $A$, then $I$ is invariant under $d_{A,A}$ because
$\Phi(d_{A,A})=[l_A,l_A]$ is contained in the Lie multiplication
algebra of $A$. Hence $\bigl(\frs\otimes I\bigr)\oplus d_{I,A}$ is
an ideal of $\frg$, and it follows that $I=A$. Hence $A$ has no
proper ideals, so it is either simple or $\dim A=1$ and $A^2=0$. In
the latter case $\Phi(\frd)=\Phi(d_{A,A})=[l_A,l_A]$ would be $0$,
and $\frg=\frs\otimes A$ would be a trivial Lie superalgebra
($[\frg,\frg]=0$), a contradiction to the simplicity of $\frg$.

Conversely, if condtions (i)--(iii) are satisfied, $A$ is unital
\cite{Shestakov}, and hence $\frs\,(\simeq \frs\otimes 1)$ is a
subalgebra of $\frg$. If $\fra$ is an ideal of $\frg$, the
invariance of $\fra$ under the adjoint action of $\frs$ shows that
$\fra=\bigl(\frs\otimes I\bigr)\oplus \fre$ for an ideal $I$ of $A$
and an ideal $\fre$ of $\frd$. Now, the simplicity of $A$ forces
that either $I=0$, but then $\fre\subseteq \ker\Phi=0$ and $\fra=0$,
or $I=A$ and then $\frs\otimes A$ is contained in $\fra$, so
$d_{A,A}=\frd$ is contained in $\fra$ too and $\fra=\frg$. Hence
$\frg$ is simple.
\end{proof}

\bigskip

Shestakov's classification \cite{Shestakov} of the simple
alternative superalgebras over $k$ (characteristic $3$) shows that
any central simple commutative alternative superalgebra is, up to
isomorphism, either:
\begin{enumerate}
\item[(i)] the ground field $k$,
\item[(ii)] the three dimensional composition superalgebra $B(1,2)$,
with even part $B(1,2)\subo=k1$,  and odd part
$B(1,2)\subuno=ku+kv$, with $1$ the unity element and $u^2=0=v^2$,
$uv=-vu=1$. This is the Jordan superalgebra of a superform on a
vector odd space of dimension $2$,
\item[(iii)] an algebra $B=B(\Gamma,D,0)$, where $\Gamma$ is a commutative
associative algebra, $D\in\der\,\Gamma$ is a derivation such that
$\Gamma$ has no proper ideal invariant under $D$, $B\subo=\Gamma$,
$B\subuno=\Gamma u$ (a copy of $\Gamma$) and the multiplication is
given by:
\begin{itemize}
\item[$\bullet$] the multiplication in $\Gamma$,
\item[$\bullet$] $a(bu)=(ab)u=(au)b$ for any $a,b\in\Gamma$,
\item[$\bullet$] $(au)(bu)=aD(b)-D(a)b$, for any $a,b\in\Gamma$.
\end{itemize}
\end{enumerate}

Given a form $\frs$ of $\frpsl_3(k)$ and the commutative alternative
superalgebra $A=B(1,2)$, its Lie superalgebra of derivations is
naturally isomorphic to $\frsl(A\subuno)\simeq\frsl_2(k)$, and the
simple Lie superalgebra $\frg=\bigl(\frs\otimes
A\bigr)\oplus\frd_{A,A}$ in Theorem \ref{th:supersAd} has even and
odd parts given by:
\[
\begin{split}
\frg\subo&=\bigl(\frs\otimes A\subo\bigr)\oplus(\frd_{A,A})\subo
  \simeq \frs\oplus\frsl_2(k),\\
\frg\subuno&=\bigl(\frs\otimes
A\subuno\bigr)\oplus(\frd_{A,A})\subuno=\frs\otimes A\subuno.
\end{split}
\]
For an algebraically closed field $k$, this coincides with the Lie
superalgebra that appears in \cite[Theorem
4.22(i)]{ElduqueNewSimple3}, and also with the derived subalgebra of
the Lie superalgebra $\frg(S_{2},S_{1,2})$ in \cite{CE1,CE2}.

\medskip

Also, assuming that $k$ is algebraically closed, according to
\cite{Block} or \cite{Yuan} any finite dimensional commutative
associative algebra $\Gamma$ over $k$ endowed with a derivation $D$
satisfying that $\Gamma$ is $D$-simple (that is, there is no proper
ideal invariant under $D$) is isomorphic to a truncated polynomial
algebra $k[t_1,\ldots,t_n: t_i^3=0,\ i=1,\ldots,n]$ which, in turn,
is isomorphic to the divided power algebra $\cO(1;n)$, which is the
$k$-algebra spanned by the symbols $t^{(r)}$, $0\leq r<3^n-1$, with
$t^{(0)}=1$ and multiplication given by
$t^{(r)}t^{(s)}=\binom{r+s}{r}t^{(r+s)}$. The isomorphism takes
$t_i$ to $t^{(3^{i-1})}$, $i=1,\ldots,n$. The simplest
$D\in\der\,\cO(1;n)$ for which $\cO(1;n)$ is $D$-simple is the
derivation given by $D:t^{(r)}\mapsto t^{(r-1)}$ for any $r$.



\section{Bouarroudj-Leites superalgebras}

Recently, S.~Bouarroudj and D.~Leites \cite{BL} have constructed an
interesting family of finite dimensional simple Lie superalgebras in
characteristic $3$ by means of the so called
Cartan-Tanaka-Shchepochkina prolongs. These superalgebras are
denoted by $\bj$, of dimension $24$, and $\Bj(1;N\vert 7)$ ($N$ an
arbitrary natural number), of dimension $2^4\times 3^N$. All these
algebras are consistently $\bZ$-graded: $\frg=\oplus_{i=-2}^{2\cdot
3^N-1}\frg_i$ and $\frg\subo$ (respectively $\frg\subuno$) is the
sum of the even (resp. odd) homogeneous components. Besides,
$\frg_0$ is the direct sum of a one dimensional center and an ideal
isomorphic to $\frpsl_3(k)$, $\dim\frg_{-2}=1$, $\frg_1$ is an
adjoint module for the ideal isomorphic to $\frpsl_3(k)$ in
$\frg_0$. The positive homogeneous components are all either a
trivial one dimensional module or an adjoint module for
$\frpsl_3(k)$, or a direct sum of both. Therefore, these Lie
superalgebras fit in the setting of the previous section.

The Lie superalgebra $\bj$ satisfies that its even part is
isomorphic to $\frsl_2(k)\oplus\frpsl_3(k)$, while its odd part is,
as a module for the even part, the tensor product of the natural two
dimensional module for $\frsl_2(k)$ and the adjoint module for
$\frpsl_3(k)$. Therefore, it coincides with the Lie superalgebra
obtained in the previous section for $A=B(1,2)$, which appeared
first in \cite{ElduqueNewSimple3}.

Over an algebraically closed field $k$ of characteristic $3$, by
dimension count, the Lie superalgebra $\Bj(1;N\vert 7)$ must be
necessarily isomorphic to a Lie superalgebra as in Theorem
\ref{th:supersAd} for $\Gamma=\cO(1;N)$ and a suitable derivation
$D$.

Let $\Gamma=\cO(1;N)$ and let $D$ be the derivation
$D:t^{(r)}\mapsto t^{(r-1)}$ for any $r$. Then the commutative
alternative superalgebra $B=B(\Gamma,D,0)=\Gamma\oplus\Gamma u$ is
consistently $\bZ$-graded with $\deg t^{(r)}=2r$, $\deg u=-1$. In
this way, $B=\oplus_{i=-1}^{2(3^N-1)}B_i$, and $\dim B_i=1$ for any
$i=-1,\ldots, 2(3^N-1)$. Then $\frd=d_{B,B}$ becomes a $\bZ$-graded
Lie superalgebra too with $\deg D=-2$. For any $a,b\in\Gamma$,
$d_{a,b}=[l_a,l_b]=0$, so $\frd\subo=[l_{B\subuno},l_{B\subuno}]$.

But for any $a,b,c\in\Gamma$:
\[
\begin{split}
[l_{au},l_{bu}](c)&=(au)(bcu)+(bu)(acu)\\
  &=D(a)bc-aD(b)c-abD(c)+D(b)ac-bD(a)c-baD(c)\\
  &=-2abD(c)=abD(c),\\[4pt]
[l_{au},l_{bu}](cu)
  &=(au)\bigl(D(b)c-bD(c)\bigr)+(bu)\bigl(D(a)c-aD(c)\bigr)\\
  &=\bigl(aD(b)c-abD(c)+D(a)bc-abD(c)\bigr)u\\
  &=D(abc)u.
\end{split}
\]
Thus $\frd\subo=\espan{d_x:x\in\Gamma}$, with $d_x\vert_\Gamma=xD$,
$d_x(yu)=D(xy)u$ for any $x,y\in \Gamma$. The degree of $d_x$ is
$\deg x -2$.

Also, $\frd\subuno=[l_{B\subo},l_{B\subuno}]$, and for any
$a,b,c\in\Gamma$:
\[
\begin{split}
[l_a,l_{bu}](c)&=a((bc)u)-(bu)(ac)=0\\[2pt]
[l_a,l_{bu}](cu)&=a\bigl(D(b)c-bD(c)\bigr)-(bu)(acu)\\
   &=aD(b)c-abD(c)-D(b)ac+bD(ac)\\
   &=(D(a)b)c.
\end{split}
\]
Thus $\frd\subuno=\espan{\delta_x:x\in\Gamma}$, with
$\delta_x\vert_\Gamma=0$, $\delta_x(yu)=xy$ for any $x,y\in\Gamma$.
The degree of $\delta_x$ is $\deg x+1$.

Hence $\frd$ has dimension $2\times 3^N$, and then the simple Lie
superalgebra $\bigl(\frpsl_3(k)\otimes B\bigr)\oplus \frd$, with
$B=B(\Gamma,D,0)$ has dimension $7\times 2\times 3^N+2\times
3^N=2^4\times 3^N$. Moreover, $\frd$ is consistently $\bZ$-graded:
\[
\frd=\frd_{-2}\oplus\frd_0\oplus\frd_1\oplus\frd_2\oplus
  \cdots\oplus\frd_{2\cdot 3^N-4}\oplus\frd_{2\cdot
  3^N-3}\oplus\frd_{2\cdot 3^N-1}.
\]
Observe that $\frd_{-1}=0=\frd_{2\cdot 3^N-2}$, while $\frd_i$ has
dimension $1$ for any other $i$ with $-2\leq i\leq 2\cdot 3^N-1$.

Now, with $\deg g=0$ for any $g\in\frpsl_3(k)$, the simple Lie
superalgebra $\frg=\bigl(\frpsl_3(k)\otimes B\bigr)\oplus\frd$ is
consistently $\bZ$-graded too:
\[
\frg=\oplus_{i=-2}^{2\cdot 3^N-1}\frg_i,
\]
with $\frg_{-2}=\frd_{-2}=kD$, $\frg_{-1}=\frpsl_3(k)\otimes
B_{-1}=\frpsl_3(k)\otimes u$, $\frg_0=\bigl(\frpsl_3(k)\otimes
B_0\bigr)\oplus\frd_0=\bigl(\frpsl_3(k)\otimes 1\bigr)\oplus
k(t^{(1)}D)$, and each $\frg_i$, $0<i<2\cdot 3^N-2$ is the direct
sum of $\frpsl_3(k)\otimes B_i$ and $\frd_i$, that is, as a a module
for $\frpsl_3(k)$, it is the direct sum of a copy of the adjoint
module and a copy of the one dimensional trivial module. Finally,
$\frg_{2\cdot 3^N-2}$ is just $\frpsl_3(k)\otimes B_{2\cdot
3^N-2}=\frpsl_3(k)\otimes t^{(3^N-1)}$ (just a copy of the adjoint
module), and $\frg_{2\cdot 3^N-1}$ is just $\frd_{2\cdot
3^N-1}=k\delta_{t^{(3^N-1)}}$ (a copy of the trivial module). This
is exactly the way $\Bj(1;N\vert 7)$ is graded, and this is no
coincidence:

\begin{theorem}
Let $k$ be an algebraically closed field of characteristic $3$, let
$N$ be a natural number, and let $\Gamma$ be the algebra of divided
powers $\cO(1;N)$. Consider the derivation $D$ of $\Gamma$ given by
$D(t^{(r)})=t^{(r-1)}$ for any $r$ and the simple commutative
alternative superalgebra $B=B(\Gamma,D,0)$. Then the simple Lie
superalgebra $\frg=\bigl(\frpsl_3(k)\otimes B\big)\oplus d_{B,B}$ in
Theorem \ref{th:supersAd} is isomorphic to the Lie superalgebra of
Bouarroudj and Leites $\Bj(1;N\vert 7)$.
\end{theorem}
\begin{proof}
Both $\frg$ and $\Bj(1;N\vert 7)$ share the same negative part
$\frg_{-}=\frg_{-2}\oplus\frg_{-1}$, and hence both of them embed in
the universal graded Lie algebra $U(\frg_{-})$ (see \cite{QG}),
which is contained in the Lie algebra of special derivations of the
tensor product of the divided power algebra $\cO(1)$ and the
Grassmann superalgebra $\Lambda(7)$ on a vector space of dimension
$7$ ($\cO(1)$ is the span of $t^{(r)}$ for any $r\geq 0$, with
$t^{(r)}t^{(s)}=\binom{r+s}{r}t^{(r+s)}$). Consider both $\frg$ and
$\Bj(1;N\vert 7)$ as subalgebras of $U=U(\frg_{-})$. Actually,
Bouarroudj and Leites consider an infinite dimensional Lie
superalgebra $\Bj(1\vert 7)$ which is contained in $U=U(\frg_{-})$.
The superalgebra $\Bj(1;N\vert 7)$ is just the intersection of
$\Bj(1\vert 7)$ with the Lie superalgebra of special derivations of
$\cO(1;N)\otimes\Lambda(7)$.

Since $\frg_{-2}=[\frg_{-1},\frg_{-1}]$,  the action of $\frg_0$ on
$\frg_{-2}$ is determined by its action on $\frg_{-1}$. By
transitivity, $U_0$ embeds in
$\End(\frg_{-1})\simeq\End(\frpsl_3(k))$, and both $\frg_0$ and
$\Bj(1;N\vert 7)_0$ act on $\frpsl_3(k)$ in the same way (the
adjoint action of $\frpsl_3(k)$ and the one dimensional center
acting as a nonzero scalar). Hence, as homogeneous subalgebras of
$U$, $\frg_0=\Bj(1;N\vert 7)_0=\Bj(1\vert 7)_0$. Write
$\Bj=\Bj(1\vert 7)$. For any $i>0$ $\Bj_i$ is defined recursively as
\[
\{x\in U_i: [x,\frg_{-2}]\subseteq \Bj_{i-2},\
[x,\frg_{-1}]\subseteq \Bj_{i-1}\},
\]
so it follows that $\frg_i\subseteq \Bj_i$. But by dimension count,
it follows that $\frg_i=\Bj(1;N\vert 7)_i=\Bj(1\vert 7)_i$ for any
$0<i<2\cdot 3^N-2$, while both $\frg_{2\cdot 3^N-2}$ and
$\Bj(1;N\vert 7)_{2\cdot 3^N-2}$ are the unique copy of the adjoint
module for $\frpsl_3(k)\subseteq U_0$ in $\Bj(1\vert 7)_{2\cdot
3^N-2}$, and both $\frg_{2\cdot 3^N-1}$ and $\Bj(1;N\vert 7)_{2\cdot
3^N-1}$ are the unique copy of the one-dimensional trivial module
for $\frpsl_3(K)$ in $\Bj(1\vert 7)_{2\cdot 3^N-1}$. Therefore,
$\frg=\Bj(1;N\vert 7)$ (as homogeneous subalgebras of $U$).
\end{proof}


\def\cprime{$'$}
\providecommand{\bysame}{\leavevmode\hbox
to3em{\hrulefill}\thinspace}
\providecommand{\MR}{\relax\ifhmode\unskip\space\fi MR }
\providecommand{\MRhref}[2]{%
  \href{http://www.ams.org/mathscinet-getitem?mr=#1}{#2}
} \providecommand{\href}[2]{#2}

\end{document}